\documentclass[12pt]{article}
\usepackage{amsmath,amsthm,amssymb,amscd}
\usepackage{latexsym}
\newtheorem{thm}{Theorem}[section]

 \newtheorem{prop}[thm]{Proposition}
 \theoremstyle{definition}
 
 \theoremstyle{remark}

 \numberwithin{equation}{section}

\title
{Sasaki manifolds with positive transverse orthogonal bisectional
curvature}

\author{ Hong Huang}
\date{}
\begin{document}
\maketitle
\begin{abstract}
 In this short note we show the following result: Let $(M^{2n+1},g)$ ($n \geq 2$) be a compact Sasaki manifold  with positive transverse orthogonal bisectional
curvature. Then $\pi_1(M)$ is finite, and the universal cover of
$(M^{2n+1},g)$ is isomorphic to a weighted Sasaki sphere. We also
get some results in the case of nonnegative transverse orthogonal
bisectional curvature under some additional conditions. This extends
recent work of He and Sun. The proof uses Sasaki-Ricci flow.


\end{abstract}
\maketitle


\section {Introduction}

In [HS1] and [HS2] He and Sun classified compact Sasaki manifolds
with nonnegative transverse bisectional curvature  using
Sasaki-Ricci flow. In this note we try to extend their results to
the case of nonnegative transverse orthogonal bisectional curvature.
First we have

 \begin{thm} \label{thm 1.1}
  Let $(M^{2n+1},g)$ ($n \geq 2$) be a compact Sasaki manifold  with positive transverse orthogonal bisectional
curvature. Then $\pi_1(M)$ is finite, and the universal cover of
$(M^{2n+1},g)$ is isomorphic to a weighted Sasaki sphere.
\end{thm}

\noindent This extends [HS1, Theorem 1.1], and is a Sasaki analogue of Gu and
Zhang [GZ, Corollary 3.2].  Note that a Sasaki manifold
$(M^{2n+1},g)$ ($n \geq 2$) has positive (resp. nonnegative)
transverse orthogonal bisectional curvature if the transverse
K$\ddot{a}$hler metric has positive (resp. nonnegative) orthogonal
bisectional curvature.

The following result extends [HS2, Theorem 2], and  is a Sasaki
analogue of  [GZ, Proposition 3.3].

\begin{thm} \label{thm 1.2}
  Let $(M^{2n+1},g)$ ($n\geq 2$) be a compact, simply connected  Sasaki manifold  with nonnegative transverse orthogonal bisectional
curvature. Suppose $b^{1,1}_B(M):=dim
H^{1,1}_B(M)=1$. Then either $M$ is isomorphic
to a weighted Sasaki sphere, or $M$ is a principal $S^1$-bundle over
an irreducible compact Hermitian symmetric space.
\end{thm}

We also have a variant of Theorem 1.2.

\begin{thm} \label{thm 1.3}
  Let $(M^{2n+1},g)$ ($n\geq 2$) be a compact, locally transversely irreducible Sasaki manifold  with nonnegative transverse orthogonal bisectional
curvature. Then either $\pi_1(M)$ is finite and the universal cover of
$(M^{2n+1},g)$ is isomorphic to a weighted Sasaki sphere, or $M$ is locally transversely symmetric.
\end{thm}

\noindent Recall  that a Sasaki manifold is locally transversely irreducible (resp. symmetric) if the transverse K$\ddot{a}$hler metric is locally irreducible (resp. symmetric). (Compare for example [HS2] and Takahashi [T].)

  In Section 2 we study Sasaki-Ricci flow with initial data a compact Sasaki manifolds with positive (nonnegative) transverse orthogonal bisectional
curvature, then in Section 3 we prove the theorems.  Since
the argument  (mainly using maximum principle) follows closely that in Chen [Ch], [GZ], [HS1] and
[HS2], the presentation will be brief.  For some notions  not
defined here see for example [HS1] and [HS2].

\section{ Sasaki-Ricci flow}

Let $(M^{2n+1}, g_0)$ be a compact Sasaki manifold. Recall [SWZ] the Sasaki-Ricci flow
\begin{equation*}
\frac{\partial  g^T}{\partial t}=-Ric^T, \hspace*{8mm}
g^T(0)=g_0^T.
\end{equation*}

 \begin{prop} \label{prop 2.1} The nonnegativity  (positivity) of transverse orthogonal bisectional curvature is preserved along the Sasaki-Ricci flow.
\end{prop}
{\bf Proof}\ \  The proof is similar to that of Proposition 2.1 in
[GZ] by using Hamilton's maximum principle and a second variation
argument originating from Mok [M]. \hfill{$\Box$

  \begin{prop} \label{prop 2.2}  \ \   Let  $(M^{2n+1},g)$ ($n\geq 2$)  be a  compact  Sasaki manifold  with nonnegative transverse orthogonal bisectional
curvature. Then any  real basic harmonic $(1,1)$-form on $M$ is
transversely parallel.  Moreover, if $b^{1,1}_B(M)=1$, then the basic first Chern class $c^B_1(M)>0$.
 \end{prop}

 {\bf Proof}\ \  The proof is similar to that of Theorem 3.1 in [GZ].  The first  conclusion is proved by using the transverse Bochner formula,
 and the second  by  using Sasaki-Ricci flow  and Proposition 2.1. \hfill{$\Box$

\vspace*{0.4cm}

\begin{prop} \label{prop 2.3} Let $(M^{2n+1},g)$ ($n\geq 2$) be a compact Sasaki manifold with
nonnegative transverse orthogonal bisectional curvature. Then
$(M^{2n+1},g)$ is transversely reducible if and only if
 $b^{1,1}_B(M)>1$.
\end{prop}

{\bf Proof}\ \  With the help of
Proposition 2.2, the proof is similar to that of Lemma 2.5 in [HS2].
 \hfill{$\Box$

\vspace*{0.4cm}

Now if $(M^{2n+1},g_0)$  ($n\geq 2$) is a  compact  Sasaki manifold
with positive transverse orthogonal bisectional curvature, then it
is locally transversely irreducible,  hence transversely
irreducible. By Proposition 2.3, $b^{1,1}_B(M)=1$. (Note that here
$b^{1,1}_B(M)=1$ can also be derived directly from the local
transverse irreducibility. Compare Theorem 3.1 (ii) in [GZ].) Then
by Proposition 2.2,  $c^B_1(M)>0$.  After a homothetic
transformation of the initial metric, we consider the normalized
Sasaki-Ricci flow

\begin{equation*}
\frac{\partial  g^T}{\partial t}=-Ric^T+g^T, \hspace*{8mm}
g^T(0)=g_0^T,
\end{equation*}
which exists on the time interval $[0,\infty)$.

\begin{prop} \label{prop 2.4}  Let $(M^{2n+1},g_0)$  ($n\geq 2$) be a  compact  Sasaki manifold  with positive transverse orthogonal bisectional
curvature.  Then along the normalized Sasaki-Ricci flow with initial
data $(M^{2n+1},g_0)$, the infimum of the transverse holomorphic
sectional curvature, if non-positive initially, will approach a
nonnegative number
as $t\rightarrow \infty$.
\end{prop}

{\bf Proof}\ \  The proof is similar to that of Theorem 1.4 in [Ch]. It is again an application of maximum principle.  \hfill{$\Box$

\begin{prop} \label{prop 2.5}  \ \  Let $(M^{2n+1},g_0)$  ($n\geq 2$) be a  compact  Sasaki manifold  with positive transverse orthogonal bisectional
curvature.  Suppose that (after  some time) $Ric^T \geq c g^T$ (for some positive constant
$c$) along the normalized Sasaki-Ricci
flow with initial data $(M^{2n+1},g_0)$. Then along the flow the
infimum of the transverse holomorphic  sectional curvature, if
non-positive initially, will become positive in finite time.
\end{prop}

{\bf Proof}\ \  The proof is similar to that of Theorem 1.5 in [Ch].  \hfill{$\Box$

\section{ Proof of Theorems}

{\bf Proof  of Theorem 1.1} \ \  Let $(M^{2n+1},g_0)$ ($n \geq 2$)
be a compact Sasaki manifold  with positive transverse orthogonal
bisectional curvature.  Then as observed above, $b^{1,1}_B(M)=1$ and
$c_1^B(M)>0$. Now we evolve the metric $g_0$ (after a homothetic
transformation if necessary) by the normalized Sasaki-Ricci flow.
Using result of Collins [Co, Theorem 1.3] and He [He, Theorem 7.1]
and Propositions 2.1, 2.4, we know that as $t \rightarrow \infty$,
$(M, g^T(t))$ subconverges to a Sasaki-Ricci soliton $(M, g_\infty)$
(on $M$) with positive transverse orthogonal bisectional curvature
and nonnegative transverse holomorphic sectional curvature, hence
nonnegative transverse bisectional curvature. It follows that  the
limit soliton  has nonnegative transverse Ricci curvature.

\vspace*{0.4cm}

 {\bf Claim} The limit soliton $(M, g_\infty)$ has
positive transverse Ricci curvature.

\vspace*{0.4cm}

\noindent Otherwise, by using Hamilton's strong maximum principle,
$Ric^T(g_\infty)$ would have a null eigenvector field $V$ which is
transversely parallel. Clearly $J_\infty V$ is also a  transversely
parallel null eigenvector field  of $Ric^T(g_\infty)$. (Here,
$J_\infty$ is the transverse complex structure of $(M, g_\infty)$.)
It follows that $(M, g_\infty)$ is locally transversely reducible. On
the other hand, $(M, g_\infty)$ has positive transverse orthogonal
bisectional curvature, and is locally transversely irreducible.
Contradiction.

\vspace*{0.4cm}

Then by Proposition 2.5,  $g^T(t)$  has positive transverse
holomorphic sectional curvature when $t
> t_0$  (for some finite $t_0$). Since $g^T(t)$  also has positive transverse orthogonal bisectional curvature, it follows that $g^T(t)$  has positive transverse
bisectional curvature when $t
> t_0$.  (Compare the proof of  [Ch, Theorem 1.8].) Using a homothetic transformation we
see  that $M$ admits a metric with $Ric> 0$, and by Myers' theorem
$\pi_1(M)$ is finite. Now the remaining conclusion in the theorem
follows from [HS1, Theorem 1.1]. \hfill{$\Box$

\vspace*{0.4cm}

{\bf Proof  of Theorem 1.2} \ \   Let $(M^{2n+1},g_0)$ ($n \geq 2$)
be a compact, simply connected Sasaki manifold with nonnegative
transverse orthogonal bisectional curvature  and with
$b^{1,1}_B(M)=1$. By Proposition 2.3, $(M^{2n+1},g_0)$ is transversely
irreducible. If $(M, g_0)$ is locally transversely symmetric, then
since we assume $M$ is simply connected,  by [T, Theorems
6.2, 6.1], $M$ is a principal $S^1$-bundle over a Hermitian symmetric
space.

\noindent Now suppose $(M, g_0)$ is not locally transversely
symmetric. We evolve the metric $g_0$ by the Sasaki-Ricci flow. Then
as in the proof of [GZ, Proposition 3.3], using a strong maximum
principle argument originating from Brendle and Schoen [BS],
Proposition 3.3 of [HS2], the assumptions on $(M, g_0)$  and
Berger's holonomy theorem, we see that the transverse orthogonal
bisectional curvature of $g(t)$ is positive for any $t\in
(0,\delta)$ for some $\delta>0$. Then the desired result follows
from Theorem 1.1.
  \hfill{$\Box$

 \vspace*{0.4cm}

{\bf Proof  of Theorem 1.3} \ \ The proof is similar to that of Theorem 1.2.  \hfill{$\Box$

\vspace*{0.4cm}

{\bf Acknowledgement}   I'm partially supported by NSFC no.11171025.

\hspace *{0.4cm}

\bibliographystyle{amsplain}

{\bf Reference}

\hspace *{0.4cm}

\bibliography{1}[BS] S. Brendle, R. Schoen, Classification of manifolds with weakly $1/4$-pinched curvature, Acta Math. 200 (2008), 1-13.

\bibliography{2}[Ch] X.X. Chen, On K$\ddot{a}$hler  manifolds  with positive orthogonal bisectional curvature, Advances in Math. 215 (2007), 427-445.

\bibliography{3}[Co] T. Collins, The transverse entropy functional
and the Sasaki-Ricci flow, Trans. Amer. Math. Soc. 365 (2013),
1277-1303.

\bibliography{4}[GZ] H. Gu, Z. Zhang, An extension of Mok's theorem on the generalized Frankel conjecture, Sci. China Math.  53 (2010), 1-12.

\bibliography{5}[He] W. He,  The Sasaki-Ricci flow and compact
Sasaki manifolds of positive transverse holomorphic bisectional
curvature, to appear in J. Geom. Anal.

\bibliography{6}[HS1] W. He, S. Sun,  Frankel conjecture and Sasaki Geometry, arXiv:1202.2589.

\bibliography{7}[HS2] W. He, S. Sun, The generalized Frankel conjecture in Sasaki
geometry, arXiv:1209.4026.

\bibliography{8}[M] N. Mok, The uniformization theorem for compact K$\ddot{a}$hler manifolds of
nonnegative holomorphic bisectional curvature, J. Diff. Geom. 27 (1988), 179-214.

\bibliography{9}[SWZ] K. Smoczyk, G-F. Wang, Y-B. Zhang. The Sasaki-Ricci flow, Internat. J.
Math. 21 (2010), no. 7, 951-969.

\bibliography{10}[T] T. Takahashi, Sasakian $\phi$-symmetric spaces,
Tohoku Math. J. 29 (1977), 91-113.

\vspace *{0.4cm}

School of Mathematical Sciences, Key Laboratory of Mathematics and Complex Systems,

Beijing Normal University, Beijing 100875, P.R. China

 E-mail address: hhuang@bnu.edu.cn

\end{document}